\documentclass[12pt]{article}
\usepackage{amsmath,amssymb,amsfonts}
\newcommand{\R}{\mathbf{R}}

\newcommand{\T}{\mathrm}

\newcommand{\fn}{\frak{n}}

\newcommand{\fm}{\frak{m}}
\newcommand{\id}{\mbox{id}}
\newcommand{\fd}{\mbox{fd}}
\newcommand{\pd}{\mbox{pd}}

\newcommand{\Gid}{\mbox{Gid}}
\newcommand{\Gfd}{\mbox{Gfd}}
\newcommand{\Gpd}{\mbox{Gpd}}

\newcommand{\para}{\paragraph}
\renewcommand{\H}{\mbox{H}}

\newcommand{\E}{\mbox{E}}
\newcommand{\K}{\mbox{K}}

\newcommand{\Hom}{\mbox{Hom}}

\renewcommand{\sup}{\mbox{sup}}
\renewcommand{\inf}{\mbox{inf}}

\begin{document}

\vspace{.5in}

\begin{center}
\LARGE{\bf Gorenstein Dimensions under~ Base~ Change}
\end{center}
\vspace{.2in}

\begin{center}
\large{\bf Leila Khatami}\footnote{lkhatami@ipm.ir} and \large{\bf
Siamak Yassemi\footnote{yassemi@ipm.ir}}\\
Department of Mathematics, University of Tehran\\and\\ Institute
for Studies in Theoretical Physics and Mathematics.
\end{center}

\vspace{.3in}

\noindent {\bf Abstract}. The so-called 'change-of-ring' results
are well-known expressions which present several connections
between projective, injective and flat dimensions over the various
base rings. In this note we extend these results to the Gorenstein
dimensions over Cohen-Macaulay local rings. \vspace{.2in}

\baselineskip=22pt

\noindent {\large \bf Introduction}

\vspace{.2in}

In 1967 Auslander [{\bf 1}] introduced the {\em Gorenstein
dimension} of a finitely generated module over a commutative
noetherian ring. This provides a characterization of Gorenstein
local rings analogous to the well--known Auslander-Buchsbaum-Serre
characterization of regular local rings. But since the Gorenstein
dimension is only defined for finitely generated modules, the
analogy is not complete. In the 1990s Enochs and Jenda in [{\bf
10}] introduced extensions of Auslander's Gorenstein dimension,
the so--called {\em Gorenstein projective} and {\em Gorenstein
flat dimensions}, and the dual notion, the {\em Gorenstein
injective dimension} and they get good results when the base ring
is Gorenstein.

Using Foxby equivalence, a nice theory for Gorenstein projective,
flat and injective dimensions over Cohen--Macaulay local rings was
given in [{\bf 13}] and [{\bf 15}].

In this note we establish connections between Gorenstein
projective, Gorenstein flat  and Gorenstein injective dimensions
of complexes over the various base rings. The note is based on the
paper " L.Khatami and S.Yassemi, {\it Gorenstein injective and
Gorenstein flat dimensions under base change}, To appear in Comm.
Algebra". Section 5 and also statements about Gorenstein
projective dimension have been added to the main paper later.

In section 1 we give the fundamental definitions and results of
hyperhomological algebra.  Details can be found in [{\bf 14}] and
this is the main reference of section 1.

In section 2 the Auslander and Bass classes over a Cohen--Macaulay
local ring with a dualizing complex are introduced. It is
well-known (cf. {\bf [4]}) that over such a ring a homologically
bounded complex has finite Gorenstein projective/flat (res.
Gorenstein injective) dimension if and only if it is in the
Auslander (res. Bass) class. Then in this section we prove some
change of ring results about Auslander and Bass classes which will
be used in our main results. For example we prove that:

Let $R$ and $S$ be $Q$--algebras such that $R$ is a Cohen-Macaulay
local ring with a dualizing complex. If $X$ and $Y$ are
homologically bounded $(R,S)$- bicomplexes and $F$ and $I$ two
$S$-complexes of finite flat and injective dimension,
respectively, then the following hold:\vspace {.1in}

(i) If $X$ belongs to the Auslander class of $R$ then $\mbox {\bf
R}\Hom_S (X,I)$ belongs to the Bass class of $R$ and $X
\otimes^{\mbox {\bf L}}_S F$ belongs to the Auslander class of
$R$.

(ii) If $Y$ belongs to the Bass class of $R$ then $\mbox {\bf
R}\Hom_S (Y,I)$ belongs to the Auslander class of $R$ and $ Y
\otimes^{\mbox {\bf L}}_S F$ belongs to the Bass class of $R$.
\vspace {.2in}

The main results of this note are proved in sections 3, 4 and 5.

In section 3 we prove the following result that is a
generalization of [{\bf 5}; 6.4.13] and [{\bf 17}; (1.4)-(1.5)].

Let $R$ and $S$ be $Q$--algebras such that $R$ is a
Cohen--Macaulay local ring with a dualizing module.If $X \in {\cal
C} (R,S)$ ($X$ is a homologically bounded complexes of
$(R,S)$--bimodules) and $Y \in {\cal C}_{(\square)} (S)$, then
\begin{itemize}
\item[(i)] $\Gpd_R (X\otimes_S^{\mathbf{L}} Y)\leq \pd_SY+\Gpd_R X$.
\item[(ii)] $\Gfd_R (X\otimes_S^{\mathbf{L}} Y)\leq \fd_SY+\Gfd_RX$.
\item[(iii)] $\Gid_R (\R \T{Hom}_S(X,Y))\leq \id_SY+\Gfd_R X$.
\item[ (iv)] If $Y \in {\cal I}(S)$, then \\$\Gfd_R ( \mbox{\bf R}\Hom_S (X,Y) ) \leq \Gid_R X +
\mbox{\rm sup} \: Y$
\end{itemize}

\vspace{.1in}

In section 4 we work on finite local ring homomorphisms of finite
flat dimension. We prove the next result which can be viewed as a
generalization of the classical results for flat and injective
dimension.

Let $ \varphi : ~ ( R ,\frak m ) \rightarrow ( S , \frak n ) $ be
a finite local ring homomorphism of Cohen--Macaulay local rings
with finite flat dimension , (that is $R$ and $ S $ are both
Cohen--Macaulay rings with unique maximal ideals $\frak m$ and
$\frak n$ respectively such that $ \varphi (\frak m) \subseteq
\frak n $ and where $S$ is a finite $R$--module with finite flat
dimension over $R$.) Then the following hold for a homologically
bounded $R$--complex $X$.

\begin{verse}

(i) $\Gpd_S ( S \otimes^{\mathbf{L}}_R X ) \leq \Gpd_R X .$

(ii) $\Gfd_S ( S \otimes^{\mathbf{L}}_R X ) \leq \Gfd_R X .$

(iii) $\Gid_S ( \mbox{\bf R}\Hom_R (S,X) ) \leq \Gid_R \: X .$

\end{verse}

In section 5 we study the connections between Gorenstein
dimensions of an $S$--complex over $R$ and $S$, when $ \phi :
(R,{\frak m}) \rightarrow (S, {\frak n})$ is a {\em quasi
Gorenstein} local ring homomorhism and $R$ a Cohen-Macaulay ring
which admits a dualizing module.

\noindent {\em Convention}. Throughout this paper by a ring we
mean a commutative noetherian ring with non-zero identity.
\vspace{.3in}

\noindent {\bf \large \bf 1.  Homological Algebra}

\vspace{.2in}

This section fixes the notation and sums up a few basic results.
The main reference is {\bf [12]} but also one can consult {\bf
[4]}.

An $R$--{\em complex} $X$ is a sequence of $R$--modules $ X_{ \ell
} $ and $ R $--linear maps $ \partial_{ \ell }^X$ , ${\ell}\in\Bbb
Z$,
\[ X = \cdots \rightarrow X_{ \ell + 1 }
\stackrel{ \partial _{ \ell + 1 } ^X }{ \rightarrow } X_{ \ell }
\stackrel{ \partial_{ \ell }^X }{ \rightarrow } X_{ \ell - 1 }
\rightarrow \cdots \] such that $ \partial _{ \ell }^X \partial _{
\ell + 1 } ^X = 0 $ for all $ {\ell} \in \Bbb Z$.  $ X_{ \ell } $
and $
\partial_{ \ell } ^X $ are called the module in degree ${\ell}$ and the
${\ell}$th differential of $X$, respectively.

The {\em supremum} and  {\em infimum} of $X$ are defined as
\[ \begin{array}{rl} \mbox{sup} \: X & \hspace{-.1in} =
\mbox{sup} \: \{ {\ell} \in \Bbb Z | \mbox{H}_{ \ell } (X) \neq 0
\}, \T{and} \\ [.1in] \mbox{inf} \: X & \hspace{-.1in} =
\mbox{inf} \: \{ {\ell} \in \Bbb Z | \mbox{H}_{ \ell } (X) \neq 0
\}.
\end{array} \]

The symbol $ {\cal C } (R) $ denotes the category of
$R$--complexes and morphisms of $R$--complexes.

The full subcategories $ {\cal C}_{ \sqsubset } (R) , {\cal C}_{
\sqsupset } (R) , {\cal C}_{\square} (R) $ and $ {\cal C}_{ 0 }
(R) $ of $ {\cal C} (R) $ consist of complexes $X$ with $ X_{\ell}
= 0 $ , for respectively $ {\ell} \gg 0 , ~ {\ell} \ll 0 , ~ |
{\ell} | \gg 0 $, and $ {\ell} \neq 0 $. The full subcategories
${\cal C}_{ ( \sqsubset ) } (R) , ~ {\cal C}_{ ( \sqsupset ) } (R)
$ and $ {\cal C}_{ (\square) } (R) $ of ${\cal C}(R)$ consist of
those $X$ with $ \mbox{H} (X) $ belonging to $ {\cal C}_{
\sqsubset } (R) , ~ {\cal C}_{ \sqsupset } (R) $, and $ {\cal
C}_{\square} (R) $, respectively.

The right derived functor of the homomorphism functor of
$R$--complexes and the left derived functor of the tensor product
of $R$--complexes are denoted by $ \mbox {\bf R} \Hom_R (-,-) $
and $ - \otimes^{\mbox {\bf L}}_R - $ , respectively.

The following inequalities hold for $ X , Z \in {\cal
C}_{(\sqsupset)} (R) $ and $ Y \in {\cal C}_{(\sqsubset)} (R)
$.(cf. {\bf [12]})

\vspace{.1in}

\[ \begin{array}{l}

 \sup ( \mbox {\bf R} \Hom_R (X,Y) )  \leqslant \sup Y - \inf X
; \mbox {and} \\ [.1in]

 \inf ( X \otimes^{\mbox {\bf L}}_R Z )  \geqslant \inf X + \inf
 Z.

 \end{array} \]

 \vspace{.1in}

A complex $X \in {\cal C}_{(\square)} (R) $ is said to be of
finite {\it {projective} }(respectively, {\it {injective}} or {\it
{flat}} ){\it {dimension}} if $ X \simeq U $,where $U$ is a
complex of projective (respectively,injective or flat) modules and
$ U_{\ell} = 0 $ for $ | \ell | \gg 0 $.

The full subcategories of $ {\cal C}_{(\square)} (R) $ consisting
of complexes of finite projective,injective and flat dimension are
denoted by $ {\cal P} (R) $ , $ {\cal I} (R) $ and $ {\cal F} (R)
$ ,respectively.

If $X$ belongs to $ {\cal C}_{(\square)} (R) $ , then the
following inequalities hold when $ P \in {\cal P} (R) $ , $ I \in
{\cal I} (R) $ and $ F \in {\cal F} (R) $.(cf. {\bf [12]})

 \vspace{.1in}

 \[ \begin{array}{l}

 \inf (\mbox {\bf R} \Hom_R (P,X) ) \geqslant  \inf X - \pd _R P;
 \\ [.1in]

 \inf (\mbox {\bf R} \Hom_R (X,I) ) \geqslant  - \sup X - \id_R I
 ;\mbox {and} \\ [.1in]

 \sup (F \otimes^{\mbox {\bf L} }_R X ) \leqslant \fd_R F + \sup
X.

\end{array} \]

\vspace{.1in}

Let $ R $ and $ S $ be commutative $Q$--algebras. Then there are
the following identities of equivalence of $Q$--complexes (cf.
section 9 of [{\bf 14}])

{\bf Commutativity.} If $ X \in {\cal C} _{ ( \sqsupset ) } (R) $
and $ Y \in {\cal C} (R) $, then
\[ X \otimes^{\mathbf{L}}_R Y = Y \otimes^{\mathbf{L}}_R X . \]

{\bf Associativity.} If $ X \in {\cal C} _{ ( \sqsupset ) } (R) ,
~ Y \in {\cal C} ( R,S ) $ and $ Z \in {\cal C} _{ ( \sqsupset ) }
(S) $, then
\[ ( X \otimes^{\mathbf{L}}_R Y )  \otimes^{\mathbf{L}}_S Z = X \otimes^{\mathbf{L}}_R
( Y \otimes^{\mathbf{L}}_S Z ) . \]

{\bf Adjointness.} If $ X \in {\cal C} _{ ( \sqsupset ) } (R) , ~
Y \in {\cal C} ( R , S ) $ and $ Z \in {\cal C} _{ ( \sqsubset ) }
(S)$, then
\[ \mbox{\bf R}\Hom_R (X, \mbox{\bf R}\Hom_S (Y, Z)) =
\mbox{\bf R}\Hom_S ( X \otimes^{\mathbf{L}}_R Y , Z) . \]

{\bf Tensor Evaluation.} If $R$ is noetherian and $ X \in {\cal C}
_{ ( \sqsupset ) } ^{ (f) }(R), ~ Y \in {\cal C}_{(\square)} ( R,S
) $ and $ Z\in {\cal C} _{ ( \sqsupset ) } (S)$, then
\[ \mbox{\bf R}\Hom_R( X, Y\otimes^{\mathbf{L}}_S Z )= \mbox{\bf R}\Hom_R(X,Y) \otimes^{\mathbf{L}}_S Z , \]
provided $ X \in {\cal P }^{(f)} (R) $ or $ Z \in {\cal F } (S) $.
\vspace{.1in}

{\bf Hom Evaluation.} If $R$ is noetherian and $ X\in {\cal C} _{
( \sqsupset ) } ^{(f)} (R) , ~ Y \in {\cal C}_{(\square)} ( R,S )
$ and $ Z\in {\cal C} _{ ( \sqsubset ) } (S)$, then
\[ X\otimes^{\mathbf{L}}_R\mbox{\bf R}\Hom_S (Y, Z)= \mbox{\bf R}\Hom_S ( \mbox{\bf R}\Hom_R (X,Y), Z ) , \]
provided $ X \in {\cal P }^{(f)} (R) $ or $ Z \in {\cal I } (S) $.

(These results actually hold under less restrictive boundedness
conditions mentioned above (and in [{\bf 14}]) but we never use
this wider generality.)

\newpage
\noindent {\large \bf 2.  The Auslander and Bass classes}

\vspace{.2in}

First recall that a complex $D$ is a {\em dualizing complex} for a
local ring $R$ when $ D \in {\cal I}^{ (f) } (R) $ and $ R =
\mbox{\bf R}\Hom_R (D,D) $.  An $R$--module $K$ which is a
dualizing complex for $R$ is said to be a {\em dualizing module}.

It is known that if $R$ admits a dualizing module, it is a
Cohen--Macaulay ring, and if $R$ is a Cohen--Macaulay ring then
every dualizing complex of $R$ has only one nonzero homology
module which can be considered as a dualizing module for $R$.
\vspace{.2in}

\noindent {\em Convention}.  In the rest of this section $(R,
\frak m ,k ) $ will be a local ring with unique maximal ideal
$\frak m$, residue field $ k $, and dualizing complex $D$.
\vspace{.2in}

The {\em Auslander class}, $ {\cal A} (R) $, of $R$ is the full
subcategory of $ {\cal C} (R) $ consisting of all $R$--complexes
$X$ satisfying
\begin{verse}

(1) $X \in {\cal C} _{(\square)} (R)$ ;

(2) $D \otimes^{\mathbf{L}}_R X \in {\cal C} _{(\square)} (R)$;
and

(3) The canonical morphism  $X\rightarrow \mbox{\bf R}\Hom_R(D, D
\otimes^{\mathbf{L}}_R X )$ is an isomorphism in ${\cal C}(R)$.
\end{verse}
\vspace{.1in}

The {\em Bass class}, $ {\cal B} (R) $, of $R$ is the full
subcategory of $ {\cal C} (R) $ consists of all $R$--complexes $Y$
satisfying

\begin{verse}

(1) $Y \in {\cal C} _{(\square)} (R)$;

(2) $\mbox{\bf R}\Hom_R (D, Y) \in {\cal C} _{(\square)} (R)$; and

(3) The canonical morphism $Y\ \leftarrow
D\otimes^{\mathbf{L}}_R\mbox{\bf R}\Hom_R (D,Y)$ is an isomorphism
in ${\cal C}(R)$.

\end{verse}

The following propositions will be used to prove the main results
of this note. \vspace{.2in}

\noindent {\bf  2.1 Proposition} Let $ R $ and $ S $ be two
$Q$--algebras. Let $X,Y\in{\cal C} _{(\square)}( R,S )$,
$F\in{\cal F}(S)$, and $I\in{\cal I}(S)$. Then the following hold:

\begin{verse}

(i) If $ X \in {\cal A} (R) $ then $ \mbox{\bf R}\Hom_S( X, I )
\in {\cal B} (R) $ and $ ( X \otimes^{\mathbf{L}}_S F ) \in {\cal
A} (R)$.

(ii) If $ Y \in {\cal B} (R) $ then $ \mbox{\bf R}\Hom_S (Y, I )
\in {\cal A} (R) $ and $ ( Y \otimes^{\mathbf{L}}_S F ) \in {\cal
B} (R) $.
\end{verse}

\vspace{.1in}

\noindent {\em Proof}.  (i) Suppose that $ X \in {\cal A} (R) $.
Since $ X \in {\cal C} _{(\square)} (R,S) $ and $ I \in {\cal I}
(S) $ and $ F \in {\cal F} (S) $, we have that $ \mbox{\bf
R}\Hom_S (X, I ) \in {\cal C} _{(\square)} (R) $ and $ X
\otimes^{\mathbf{L}}_S F \in {\cal C} _{(\square)} (R) $. From the
equalities $$\mbox{\bf R}\Hom_R (D,\mbox{\bf R}\Hom_S (X, I))=
\mbox{\bf R}\Hom_S ( D \otimes^{\mathbf{L}}_R X , I) \mbox{;
and}$$ $$ D \otimes^{\mathbf{L}}_R (X \otimes^{\mathbf{L}}_S F) =
(D \otimes^{\mathbf{L}}_R X) \otimes^{\mathbf{L}}_S F $$ we get
that $\mbox{\bf R}\Hom_R (D,\mbox{\bf R}\Hom_S (X, I))$ and $ D
\otimes^{\mathbf{L}}_R (X \otimes^{\mathbf{L}}_S F) $ are both
homologically bounded since so is $ D \otimes^{\mathbf{L}}_R X $.

Finally note that
\[ \begin{array}{rl}
D\otimes^{\mathbf{L}}_R(\mbox{\bf R}\Hom_R(D, \mbox{\bf R}\Hom_S
(X,I )) &\, = D\otimes^{\mathbf{L}}_R(\mbox{\bf R}\Hom_S ( D
\otimes^{\mathbf{L}}_R X, I ))\\ &\, = \mbox{\bf R}\Hom_S(
\mbox{\bf R}\Hom_R(D,D \otimes^{\mathbf{L}}_R X),I)
\end{array}\]
which is canonically isomorphic to $ \mbox{\bf R}\Hom_S (X,I) $;
and
\[\begin{array}{rl}
\mbox{\bf R}\Hom_R(D,D \otimes^{\mathbf{L}}_R (X
\otimes^{\mathbf{L}}_S F ))&\, =\mbox{\bf
R}\Hom_R(D,(D\otimes^{\mathbf{L}}_R X)
\otimes^{\mathbf{L}}_S F)\\
&\, = \mbox{\bf R}\Hom_R(D,D \otimes^{\mathbf{L}}_R X)
\otimes^{\mathbf{L}}_S F
\end{array}\] which is canonically isomorphism with $ X \otimes^{\mathbf{L}}_S F $.

Therefore $ \mbox{\bf R}\Hom_S (X,I) \in {\cal B} (R) $ and $
X\otimes^{\mathbf{L}}_S F \in {\cal A} (R) $.

\noindent The proof of part (ii) is analogous to the proof of part
(i).\hfill$\square$

Recall that an injective $S$--module $J$ is called {\em faithfully
injective} if for all non-zero $S$--modules $M$ we have $ \Hom_S (
M,J) \neq 0 $.

A flat $S$--module $P$ is called {\em faithfully flat} if for all
nonzero $S$--modules $M$ we have $ M \otimes_S P \neq 0 $.

Since in this note all rings are assumed to be noetherian it is
clear that the module $ J = \coprod_{ \frak n \in \mbox{max}S } \:
\E_S ( S/\frak n ) $, the sum of injective envelopes of the $ S /
\frak n $, for all maximal ideals of $S$, is a faithfully
injective $S$--module and $S$ itself is a faithfully flat
$S$--module.\vspace{.2in}

\noindent {\bf 2.2 Proposition}  With the same conditions as those
of the proposition (2.1), if $J$ is a faithfully injective
$S$--module and $P$ is a faithfully flat $S$--module, then the
following hold:

\begin{verse}
(i) $X \in {\cal A}  (R)$ if and only if $ \Hom_S ( X,J ) \in
{\cal B} (R) $

(ii) $X \in {\cal A}  (R)$ if and only if $ X \otimes_S P \in
{\cal A}  (R) $

(iii) $X \in {\cal B}  (R)$ if and only if $ \Hom_S ( Y,J ) \in
{\cal A} (R) $

(iv) $X \in {\cal A}  (R)$ if and only if $ Y \otimes_S P \in
{\cal B}  (R) $.
\end{verse}
\vspace{.2in}

\noindent {\em Proof}. The ``only if'' parts are clear from
Proposition (2.1) since $ \mbox{\bf R}\Hom_S( Z, J ) $ and $ Z
\otimes^{\mathbf{L}}_S P $ are isomorphic to $ \Hom_S ( Z , J ) $
and $ Z \otimes_S P $ respectively, for all $ Z\in {\cal C} (S) $.
\vspace{.1in}

Now we prove the ``if'' part of (i). The other claims can be
proved with similar techniques.

Since $J$ is a faithfully injective $S$--module, if $ \Hom_S ( X,J
) $ is homologically bounded then so is $X$.  We also have
\[\begin{array}{rl}
\Hom_S(D\otimes^{\mathbf{L}}_R X, J ) &\, = \mbox{\bf R}\Hom_S ( D\otimes^{\mathbf{L}}_R X, J )\\
&\, = \mbox{\bf R}\Hom_R (D, \Hom_S ( X , J )) \in {\cal C}
_{(\square)} (R)
\end{array}\] then $ D \otimes^{\mathbf{L}}_R X \in {\cal C} _{(\square)} (R) $. Finally
$$\begin{array}{rl}
\Hom_S (\mbox{\bf R}\Hom_R (D,D \otimes^{\mathbf{L}}_R X ), J
)&\, = \mbox{\bf R}\Hom_S (\mbox{\bf R}\Hom_R (D,D \otimes^{\mathbf{L}}_R X ), J )\\
&\, = D\otimes^{\mathbf{L}}_R\mbox{\bf R}\Hom_S (D
\otimes^{\mathbf{L}}_R X,
J)\\ &\, = D\otimes^{\mathbf{L}}_R\mbox{\bf R}\Hom_R (D, \mbox{\bf R}\Hom_S (X, J ))\\
&\, = D\otimes^{\mathbf{L}}_R(\mbox{\bf R}\Hom_R(D, \Hom_S ( X,J
)).
\end{array}$$

\noindent Then $ \Hom_S (\mbox{\bf R}\Hom_R (D,D
\otimes^{\mathbf{L}}_R X ), J ) $ is canonically isomorphism with
$ \Hom_S ( X,J ) $ and hence $ \mbox{\bf R}\Hom_R (D,
D\otimes^{\mathbf{L}}_RX) $ is canonically isomorphism with $X$.
\hfill $ \square $ \vspace{.2in}

\noindent {\bf 2.3 Proposition} Let $ \varphi : ( R , \frak m )
\rightarrow ( S , \frak n ) $ be a finite local homomorphism of
local rings such that $S$ has finite flat dimension over $R$. If
$R$ has a dualizing complex $D$, then for $ X \in {\cal C}
_{(\square)} (R) $ the following hold:

\begin{verse}
(i) If $ X \in {\cal A}  (R) $ then $ S \otimes^{\mathbf{L}}_R X
 \in {\cal A} (S) $

(ii) If $ X \in {\cal B}  (R) $ then $\mbox{\bf R}\Hom_R (S,X )
\in {\cal B} (S) $.
\end{verse}
\vspace{.1in}

\noindent {\em Proof}. Note that the $S$--complex $ \tilde D=
\mbox{\bf R}\Hom_R (S,D) $ is a dualizing complex for $S$, cf.
[{\bf 14}; 15.28]. Also recall that since $R$ is noetherian and
$S$ has finite flat dimension over $R$ it has finite projective
dimension over $R$ too.

\noindent (i): See [{\bf 6}; 6.6].

\vspace{.1in}

\noindent (ii): It is clear that $ \mbox{\bf R}\Hom_R (S,X) \in
{\cal C} _{(\square)} (S) $.  We also have the following
equalities:
\[ \begin{array}{rl}
\mbox{\bf R}\Hom_S(\tilde D,\mbox{\bf R}\Hom_R(S,X))&\, = \mbox{\bf R}\Hom_R (\tilde D \otimes^{\mathbf{L}}_SS, X)\\
&\, =\mbox{\bf R}\Hom_R (\tilde D, X )
\\ &\, = S\otimes^{\mathbf{L}}_R \mbox{\bf R}\Hom_R (D,X )
\end{array} \] The latter is homologically bounded since $\mbox{\bf R}\Hom_R (D,X ) \in {\cal C} _{(\square)} (R) $ and $ \fd_R
S<\infty$, hence $ \mbox{\bf R}\Hom_S(\tilde D,\mbox{\bf
R}\Hom_R(S,X)) \in {\cal C} _{(\square)} (S) $. Finally note that
\[ \begin{array}{rl}
\tilde D\otimes^{\mathbf{L}}_S\mbox{\bf R}\Hom_S(\tilde
D,\mbox{\bf R}\Hom_R(S,X))&\, =\tilde
D\otimes^{\mathbf{L}}_S(\mbox{\bf R}\Hom_R(\tilde D, X))
\\ &\, = \tilde D\otimes^{\mathbf{L}}_R(\mbox{\bf R}\Hom_R(D,X))\\
&\, =\mbox{\bf R}\Hom_R(S,D\otimes^{\mathbf{L}}_R(\mbox{\bf
R}\Hom_R(D,X))).
\end{array} \] The complex $X$ represents $D\otimes^{\mathbf{L}}_R\mbox{\bf R}\Hom_R(D,X) $ canonically,
then $ \mbox{\bf R}\Hom_R (S,X) $ canonically represents $\tilde
D\otimes^{\mathbf{L}}_S\mbox{\bf R}\Hom_S(\tilde D,\mbox{\bf
R}\Hom_R(S,X))$. \hfill $ \square $ \vspace{.3in}

\noindent {\large \bf 3.  Gorenstein dimensions}

\vspace{.2in}

In this section we study the Gorenstein dimensions of complexes of
modules over a Cohen--Macaulay local ring which admits a dualizing
module.

Let $ P \in {\cal C} ^P (R) $ be homologically trivial. $P$ is
called a {\em complete projective resolution} if the complex $
\Hom_R(P,Q) $ is homologically trivial for every projective
$R$--module $Q$.

A module $M$ is said to be {\em Gorenstein projective} if there
exists a complete projective resolution $P$ with $ C_0^P \cong M
$. Observe that every projective module is obviously Gorenstein
projective.

The notation $ {\cal C} ^{ GP } (R) $ is used for the full
subcategory (of $ {\cal C} (R) $) of complexes of Gorenstein
projective modules.

The {\em Gorenstein projective dimension} of $ X \in {\cal C}_{
(\sqsupset) } (R) $, $\Gpd_R X $, is defined as
\[ \Gpd_R X  = \mbox{inf} \: \{\: \mbox{sup} \: \{\: l \in
\Bbb Z | A_{ \ell } \neq 0 \} \:| \: X \simeq A \in {\cal C}_{
\sqsupset } ^{ GP } (R) \: \} . \] (The set over which infimum is
taken is non--empty since any complex $X\in {\cal C}_{ (\sqsupset)
} (R)$ has a projective resolution that belongs to ${\cal C}_{
\sqsupset } ^{ GP } (R)$.)

In particular, if $M$ is a non--zero $R$--module then $\Gpd_RM$ is
the smallest integer $n\ge 0$, such that there is an exact
sequence
$$0\to T_n\to\cdots\to T_1\to T_0\to M\to 0$$ where each $T_i$ is
a non--zero Gorenstein projective module.

\noindent {\bf 3.1 Theorem} ([{\bf 5}], 4.4.5--4.4.16) Let $R$ be
a Cohen--Macaulay local ring which admits a dualizing module. For
a complex $ X \in {\cal C} _{ ( \sqsupset ) } (R) $ the next three
conditions are equivalent.

\begin{verse}
(i)  $ X \in {\cal A}  (R) $

(ii) $ \Gpd_R X < \infty $

(iii) $ X \in {\cal C} _{(\square)} (R) $ and $ \Gpd_R X \leq \sup
X + \dim R $.
\end{verse}

Furthermore, if $ X \in {\cal A}  (R) $, then
\[ \begin{array}{rl}
\Gpd_R X & \hspace{-.1in} = \mbox{\rm sup} \: \{ \mbox{\rm inf} \:
U  - \mbox{\rm inf} \: (\T{\mathbf R}\Hom_R(X,U)) | U \in {\cal F}
(R) \wedge U \not \simeq 0 \} \\ [.1in] & \hspace{-.1in} =
\mbox{\rm sup} \: \{ -\mbox{\rm inf} \: (\T{\mathbf R}\Hom_R(X,Q))
| Q \in {\cal C}_0^P (R) \}
\\ [.1in] & \hspace{-.1in} = \mbox{\rm sup} \: \{ \mbox{\rm inf} \:
U  - \mbox{\rm inf} \: (\T{\mathbf R}\Hom_R(X,U)) | U \in {\cal I}
(R) \wedge U \not \simeq 0 \} \\ [.1in] & \hspace{-.1in} =
\mbox{\rm sup} \: \{ -\mbox{\rm inf} \: (\T{\mathbf R}\Hom_R(X,T))
| T \in {\cal I}_0 (R) \}.
\end{array} \]

Note that it is clear from the definition of the Gorenstein
projective dimension that the inequality $ \Gpd \: X \leq \pd_R X
$holds. It is known that if $ \pd_R X < \infty $ then equality
holds. (cf.\ [{\bf 5}; (4.4.7)].)

Let $ F \in {\cal C} ^F (R) $ be homologically trivial. $F$ is
called a {\em complete flat resolution} if the complex $ J
\otimes_R F $ is homologically trivial for every injective
$R$--module $J$.

A module $M$ is said to be {\em Gorenstein flat} if there exists a
complete flat resolution $F$ with $ C_0^F \cong M $. Observe that
every flat module is obviously Gorenstein flat.

The notation $ {\cal C} ^{ GF } (R) $ is used for the full
subcategory (of $ {\cal C} (R) $) of complexes of Gorenstein flat
modules.

The {\em Gorenstein flat dimension} of $ X \in {\cal C}_{
(\sqsupset) } (R) $, $\Gfd_R X $, is defined as
\[ \Gfd_R X  = \mbox{inf} \: \{\: \mbox{sup} \: \{\: l \in
\Bbb Z | A_{ \ell } \neq 0 \} \:| \: X \simeq A \in {\cal C}_{
\sqsupset } ^{ GF } (R) \: \} . \] (The set over which infimum is
taken is non--empty since any complex $X\in {\cal C}_{ (\sqsupset)
} (R)$ has a projective resolution that belongs to ${\cal C}_{
\sqsupset } ^{ GF } (R)$.)

In particular, if $M$ is a non--zero $R$--module then $\Gfd_RM$ is
the smallest integer $n\ge 0$, such that there is an exact
sequence
$$0\to T_n\to\cdots\to T_1\to T_0\to M\to 0$$ where each $T_i$ is
a non--zero Gorenstein flat module.

\vspace{.2in}

\noindent {\bf 3.2 Theorem} ([{\bf 5}], 5.2.6--5.4.6) Let $R$ be a
Cohen--Macaulay local ring which admits a dualizing module. For a
complex $ X \in {\cal C} _{ ( \sqsupset ) } (R) $ the next three
conditions are equivalent.

\begin{verse}
(i)  $ X \in {\cal A}  (R) $

(ii) $ \Gfd_R X < \infty $

(iii) $ X \in {\cal C} _{(\square)} (R) $ and $ \Gfd_R X \leq \sup
X + \dim R $.
\end{verse}

Furthermore, if $ X \in {\cal A}  (R) $, then
\[ \begin{array}{rl}
\Gfd_R X & \hspace{-.1in} = \mbox{\rm sup} \: \{ \mbox{\rm sup} \:
( U \otimes^{\mathbf{L}}_R X ) - \mbox{\rm sup} \: U | U \in {\cal
I} (R) \wedge U \not \simeq 0 \} \\ [.1in] & \hspace{-.1in} =
\mbox{\rm sup} \: \{ \mbox{\rm sup} \: ( U \otimes^{\mathbf{L}}_R
X ) - \mbox{\rm sup} \: U | U \in {\cal F} (R) \wedge U \not
\simeq 0 \} \\ [.1in]& \hspace{-.1in} = \mbox{\rm sup} \: \{
\mbox{\rm sup} \: ( J \otimes^{\mathbf{L}}_R X ) | J \in {\cal
C}_0^I (R) \} \\ [.1in] & \hspace{-.1in} = \mbox{\rm sup} \: \{
\mbox{\rm sup} \: ( T \otimes^{\mathbf{L}}_R X ) | T \in {\cal
I}_0^f (R) \}.
\end{array} \]

Note that it is clear from the definition of the Gorenstein flat
dimension that the inequality $ \Gfd \: X \leq \fd_R X $ holds. It
is known that if $ \fd_R X < \infty $ then equality holds. (cf.\
[{\bf 5}; (5.2.9)].)

Let $ I \in {\cal C} ^I (R) $ be homologically trivial. $I$ is
called a {\em complete injective resolution} if the complex $
\Hom_R ( J , I ) $ is homologically trivial for every injective
$R$--module $J$.

A module $N$ is said to be {\em Gorenstein injective} if there
exists a complete injective resolution $I$ with $ Z_0^I \cong N $.
Observe that every injective module is obviously Gorenstein
injective.

The notation $ {\cal C} ^{ G I } (R) $ is used for the full
subcategory (of $ {\cal C} (R) $) of complexes of Gorenstein
injective modules.

The {\em Gorenstein injective dimension} of $ Y \in {\cal C} _{ (
\sqsubset ) } (R) $ , $\Gid_R Y $, is defined as
\[ \Gid_R Y = \mbox{inf} \:  \{\: \mbox{sup} \:
\{\: {\ell} \in \Bbb Z | B_{ - {\ell} } \neq 0 \: \}\: | \: Y
\simeq B \in {\cal C}_{ \sqsubset } ^{ G I } (R) \}\: \]

(The set over which infimum is taken is non--empty since any
complex $ Y \in {\cal C} _{ ( \sqsubset ) } (R) $ has an injective
resolution that belongs to ${\cal C}_{ \sqsubset } ^{ GI } (R)$.)

So in particular, if $M$ is a non--zero $R$--module then $\Gid_RM$
is the smallest integer $n\ge 0$, such that there is an exact
sequence
$$0\to M\to\H^0\to H^1\to \cdots\to H^n\to 0$$ where each $H^i$ is
a non--zero Gorenstein injective module.

\vspace{.2in}

\noindent {\bf 3.3 Theorem} ([{\bf 5}], 6.2.5) Let $R$ be a
Cohen--Macaulay local ring which admits a dualizing module. For a
complex $ Y \in {\cal C} _{ ( \sqsubset ) } (R) $ the following
conditions are equivalent.

\begin{verse}
(i) $ Y \in B (R) $

(ii) $ \Gid_R Y < \infty $

(iii) $ Y \in C _{(\square)} (R) $ and $ \Gid_R Y \leq - \inf Y +
\dim R $.
\end{verse}

Furthermore if $ Y \in B (R) $, then
\[ \begin{array}{rl}
\Gid_R Y  & \hspace{-.1in} = \mbox{\rm sup} \: \{ - \mbox{\rm sup}
\: U - \mbox{\rm inf} \:  ( \mbox{\bf R}\Hom_R (U,Y) ) | U \in
{\cal I} (R) \wedge U \not \simeq 0 \}
\\ [.1in]
& \hspace{-.1in} = \mbox{\rm sup} \: \{ - \mbox{\rm inf} \:  (
\mbox{\bf R}\Hom_R (J,Y) ) | J \in {\cal C} _0^I (R) \} .
\end{array}
\]

\vspace{.1in}

Note that from the definition of the Gorenstein injective
dimension that the inequality $ \Gid_R Y \leq \id_R Y $ holds. It
is known that over a local ring with dualizing module if $ \id_R Y
< \infty $ then equality holds. (cf.\ [{\bf 5}; (6.2.6)].)
\vspace{.2in}


In the next theorem $R$ and $S$ are assumed to be $Q$--algebras.
\vspace{.2in}

\noindent {\bf 3.4 Theorem} Let $R$ be a Cohen--Macaulay local
ring with a dualizing module $D$. If $X \in {\cal C} (R,S)$ and $Y
\in {\cal C}_{(\square)} (S)$, then
\begin{itemize}
\item[(i)] $\Gpd_R (X\otimes_S^{\mathbf{L}} Y)\leq \pd_SY+\Gpd_R X$.
\item[(ii)] $\Gfd_R (X\otimes_S^{\mathbf{L}} Y)\leq \fd_SY+\Gfd_RX$.
\item[(iii)] $\Gid_R (\R \T{Hom}_S(X,Y))\leq \id_SY+\Gfd_R X$.
\item[ (iv)] If $Y \in {\cal I}(S)$, then \\$\Gfd_R ( \mbox{\bf R}\Hom_S (X,Y) ) \leq \Gid_R X +
\mbox{\rm sup} \: Y$.
\end{itemize}

\vspace{.1in}

\noindent {\em Proof}. To prove each inequality we assume that the
right hand side terms are finite, the inequality is clear if one
of the terms is infinite.

By (3.1, 3.2, 3.3) and (2.1), finiteness of the right hand side
terms of each inequality implies the finiteness of the left hand
term side of it.

(i) We have
\[ \begin{array}{rl} \Gpd_R ( X \otimes^{\mathbf{L}}_S Y )& \, = \sup\{ - \inf( \T{\mathbf R}\Hom_R(X
\otimes^{\mathbf{L}}_S Y,T) | T \in {\cal C}_0^P (S) \}
\\ &\, = \sup\{ - \inf( \T{\mathbf R}\Hom_S(Y,\T{\mathbf R}\Hom_R(X,T))
| T \in {\cal C} _0^P (R) \}
\\ & \,\leq\sup\{ -\inf(\T{\mathbf R}\Hom_R(X,T)) + \pd_S Y | T \in
{\cal C} _0^P (R) \}
\\
&\, = \Gpd_R X + \pd _S Y .
\end{array} \]

(ii) We have
\[ \begin{array}{rl} \Gfd_R ( X \otimes^{\mathbf{L}}_S Y )& \, = \sup\{ \sup( J\otimes^{\mathbf{L}}_R( X
\otimes^{\mathbf{L}}_S Y )_R ) | J \in {\cal C}_0^I (R) \}
\\ &\, = \sup\{ \sup(( J \otimes^{\mathbf{L}}_R X) \otimes^{\mathbf{L}}_S Y )) | J \in
{\cal C} _0^I (R) \}
\\ & \,\leq\sup\{ \sup( J \otimes^{\mathbf{L}}_R X ) + \fd_S Y | J \in
{\cal C} _0^I (R) \}
\\
&\, = \Gfd_R X + \fd _S Y .
\end{array} \]

(iii) We have
\[ \begin{array}{rl} \Gid_R (\mbox{\bf R}\Hom_S (X,Y)
) &\, =\sup\{-\inf(( \mbox{\bf R}\Hom_R (J,\mbox{\bf R}\Hom_S
(X,Y)) | J \in {\cal C}_0^I (R) \}\\ &\, =\sup\{-\inf( \mbox{\bf
R}\Hom_S (
J\otimes^{\mathbf{L}}_R X,Y) ) | J \in {\cal C}_0^I (R) \}\\
&\, \le\sup\{ \sup(J \otimes^{\mathbf{L}}_R X ) + \id_S Y | J \in
{\cal C}_0^I
(R) \}\\
&\,  = \Gfd_R X + \id_S Y .
\end{array} \]

(iv) We have
\[ \begin{array}{rl}
\Gfd_R ( \mbox{\bf R}\Hom_S (X,Y )) &\, =\sup\{ \sup(
(J\otimes^{\mathbf{L}}_R\mbox{\bf R}\Hom_S (X,Y))| J \in {\cal
I}_0^f (R) \}
\\ &\, = \sup\{ \sup(\mbox{\bf R}\Hom_S ( \mbox{\bf R}\Hom_R (J,X ), Y)) | J \in
{\cal I} _0^f (R) \}
\\
&\leq\sup\{ - \inf(\mbox{\bf R}\Hom_R (J,X ) + \sup Y | J \in
{\cal I} _0^f (R) \}
\\ &\, \le\Gid_R X + \sup Y .
\end{array} \] \hfill $ \square $ \vspace{.1in}

\noindent {\bf 3.5 Corollary} Let $R$ be a Cohen-Macaulay local
ring with a dualizing module. If $\phi : R \rightarrow S$ is a
ring homomorphism and $Y \in {\cal C}_{(\square)}(S)$, then
\begin{itemize}
\item[(i)] $\Gpd_R Y\leq \pd_S Y+\Gpd_R S$.
\item[(ii)] $\Gfd_R Y\leq \fd_S Y+\Gfd_R S$.
\item[(iii)] $\Gid_R Y\leq \id_S Y+\Gfd_R S$.
\item[(iv)] If $\id_R \: Y < \infty$ , then $\Gfd_R \: Y \leq \Gid_R \:
S + \mbox{sup} \: Y$.\hfill $ \square $
\end{itemize}

\vspace{.2in}

\noindent {\bf 3.6 Corollary} With the same conditions as (3.5),
if $S$ is a Gorenstein local ring then
\[ \Gfd_R \: S \leq \Gid_R \: S
\leq \Gfd_R \: S + \mbox{\rm dim}\: S . \] In particular if $S$ is
self-injective, then $ \Gfd_R \: S = \Gid_R \: S $. \hfill $
\square $ \vspace{.2in}

\noindent {\bf 3.7 Proposition} Let $R$ and $S$ be $Q$--algebras
such that $R$ is a Cohen--Macaulay local ring with a dualizing
module. If $J$ is a faithfully injective $S$--module and $P$ is a
faithfully flat $S$--module, then the following hold for $ X \in
{\cal C} _{(\square)} ( R,S ) $.

\begin{itemize}
\item[(i)]$\Gfd_R ( X \otimes_S P ) = \Gfd_R \: X$
\item[(ii)]$ \Gid_R ( \Hom_S ( X,J ) ) = \Gfd_R \: X$
\item[(iii)]$\Gfd_R ( \Hom_S ( X,J )) < \infty$ then $\Gid_R X < \infty$.
\end{itemize}


\noindent {\em Proof}.  (i), (ii) By (2.2(i),(ii)), (3.2) and
(3.3) we have that $ \Gid_R \Hom_S ( X,J )<\infty $ if and only if
$ \Gfd_R \: X < \infty $ if and only if $ \Gfd_R ( X \otimes_S P )
< \infty $.

The equalities are clear from the proof of (3.4) , since the
inequalities there become equalities for faithfully injective $J$
and faithfully flat $P$. \vspace{.1in}

\noindent (iii) Use (3.2),(3.3) and (2.2 (iii)).\hfill $ \square $
\vspace{.1in}

\noindent {\it Remark.} It is natural to ask about the equality of
part (c). We were unable to prove that equality always holds. We
do not even know whether an $R$--module $M$ is Gorenstein
injective if and only if $\Hom_R(M,E)$ is Gorenstein flat, for
every injective $R$--module $E$. \vspace{.2in}

\noindent {\bf 3.8 Corollary}   Let $R$ be a Cohen--Macaulay local
ring with a dualizing module and $ \varphi : ~ R \rightarrow S $ a
ring homomorphism. If $ X \in {\cal C} _{(\square)} (S) $, and
$\frak p$ is a prime ideal of $S$
\begin{itemize}
\item[(i)] If $ \Gfd_R X _{\frak p} < \infty $, then
$$ \Gid_R ( \Hom_S ( X , \mbox{E}_S ( S/\frak p ))) = \Gfd_R
X_{\frak p} .$$

\item[(ii)] If $ \Gid_R X_{\frak p} < \infty $, then $$\Gfd_R (
\Hom_S ( X ,\mbox{E}_S ( S/\frak p ))) \leq \Gid_R \: X_{\frak p}
.$$
\end{itemize}
\noindent {\em Proof}.  Note that
\[ \begin{array}{rl}
\Hom_S ( X , \mbox{E}_S ( S/\frak p )) & \, = \Hom_S ( X ,
\Hom_{S_{\frak
p}} ( S_{\frak p} , \E_S ( S/{\frak p} ))\\
& \,= \Hom_{S_{\frak p}} ( X_{\frak p} , \E_S ( S/{\frak p} )) .
\end{array} \] Now use (3.7). \hfill $ \square $ \vspace{.3in}

\noindent {\large \bf 4.  Finite local ring homomorphism}

\vspace{.2in}

In this section $ \varphi : ~ ( R ,\frak m ) \rightarrow ( S ,
\frak n ) $ is a finite local ring homomorphism of Cohen--Macaulay
local rings.We also assume that $R$ has a dualizing module $D$,
and consequently $ (\mbox{\bf R}\Hom_R (S, D)) $ is a dualizing
module for $S$. \vspace{.2in}

\noindent {\bf 4.1 Theorem}  If $ \mbox{\rm fd}_R \: S < \infty $,
then the following hold for $ X \in {\cal C} _{(\square)} (R) $.

\begin{verse}

(i) $\Gpd_S ( S \otimes^{\mathbf{L}}_R X ) \leq \Gpd_R X .$

(ii) $\Gfd_S ( S \otimes^{\mathbf{L}}_R X ) \leq \Gfd_R X .$

(iii) $\Gid_S ( \mbox{\bf R}\Hom_R (S,X) ) \leq \Gid_R \: X .$

\end{verse}

\vspace{.1in}

\noindent {\em Proof}.  We prove each inequality when the right
hand side term of it is finite, otherwise there is nothing to
prove.

(i) We assume that $ \Gpd_R \: X < \infty $. Then by (2.3(i)) $
\Gpd_S ( S \otimes^{\mathbf{L}}_R X ) < \infty $.

Therefore
\[ \begin{array}{rl}
\Gpd_S ( S \otimes^{\mathbf{L}}_R X ) &\, = \sup\{ -\inf
(\mbox{\bf R}\Hom_S( S \otimes^{\mathbf{L}}_R X ,Q))| Q \in {\cal
C} _0^P (S) \} \\&\, = \sup\{ -\inf(\mbox{\bf R}\Hom_R(X,\mbox{\bf
R}\Hom_S(S,Q))) | Q\in {\cal C} _0^P (S) \}
\\&\, = \sup\{-\inf(\mbox{\bf R}\Hom_R(X,Q))| Q\in {\cal C} _0^P (S) \}
\end{array}
\]

By [{\bf 2}; (4.2(b))] $Q$ has finite projective dimension over
$R$ if it has finite projective dimension over $S$.  Then
\[ \Gpd_S ( S\otimes^{\mathbf{L}}_R X ) \leq \Gpd_R \: X . \]

(ii) We assume that $ \Gfd_R \: X < \infty $. Then by (2.3(i)) $
\Gfd_S ( S \otimes^{\mathbf{L}}_R X ) < \infty $.

Therefore
\[ \begin{array}{rl}
\Gfd_S ( S \otimes^{\mathbf{L}}_R X ) &\, = \sup\{
\sup(J\otimes^{\mathbf{L}}_S( S \otimes^{\mathbf{L}}_R X ))| J \in
{\cal C} _0^I (S) \}
\\
&\, = \sup\{ \sup( J \otimes^{\mathbf{L}}_R X ) | J \in {\cal C}
_0^I (S) \}
\end{array}
\]

By [{\bf 2}; (4.2(b))] $J$ has finite injective dimension over $R$
if it has finite injective dimension over $S$.  Then
\[ \Gfd_S ( S\otimes^{\mathbf{L}}_R X ) \leq \Gfd_R \: X . \]

(iii) We assume that $ \Gid_R \: X < \infty $.  Then by (2.3(ii))
$ \Gid_S (\mbox{\bf R}\Hom_R (S,X) ) < \infty $.  And
\[ \begin{array}{rl}
\Gid_S( \mbox{\bf R}\Hom_R (S,X) ) &\, =\sup\{ - \inf(\mbox{\bf
R}\Hom_S(J, \mbox{\bf R}\Hom_R(S,X))| J \in {\cal C}_0^I (S) \}
\\ &\,= \sup\{ - \inf(\mbox{\bf R}\Hom_R (J,X) ) | J \in
{\cal C} _0^I (S) \}
\\ &\,\leq \Gid_R \: X .  \end{array} \] \hfill $ \square $
\vspace{.2in}

Recall that for $ x \in\frak m $ the complex $ 0 \rightarrow R
\stackrel{ \cdot x }{ \rightarrow } R \rightarrow 0 $ concentrated
in degrees one and zero is called the {\em Koszul complex} of $x$
and denoted $\K(x) $.  For $ x = x_1 ,..., x_n \in\frak m $ the
Koszul complex $\K(x) $ of $ x_1 ,..., x_n $ is the complex $\K(
x_1 ) \mbox{$ \otimes_R $} \cdots \mbox{$ \otimes_R $}\K(x_n) $.
Note that $\K(x)$ is a homologically bounded $R$--complex of
finite free $R$--modules.

The sequence $ x = x_1 ,..., x_n $ of elements of $\frak m$ are
$R$--regular if and only if $\K(x) $ has zero homology modules
except at $\K(x)_0$,and when this is the case the homology module
in degree zero is $ R/(x) $. (cf. [{\bf 4}; 1.6.19]).
\vspace{.2in}

\noindent {\bf 4.2 Proposition} Let $ x = x_1 ,..., x_n \in R $ be
an $R$--regular sequence and let $ X \in {\cal C} _{(\square)} (R)
$, then the following hold

\begin{verse}

(i) $\Gpd_{ R/(x) } \: ( K(x) \otimes_R X ) \leq \Gpd_R \: X $

(ii) $\Gfd_{ R/(x) } \: ( K(x) \otimes_R X ) \leq \Gfd_R \: X $

(ii) $\Gid_{ R/(x) } \: ( \Hom _R ( K(x) , X )) \leq \Gid_R \: X
.$

\end{verse}

\vspace{.1in}

\noindent {\em Proof}.  It is clear from remark (4.2) that $ K(x)
$ is a free resolution for $ R / (x) $, then the above
inequalities are consequences of (4.1). \hfill $ \square $
\vspace{.1in}

\noindent {\it Remark.} Note that the converses inequalities of
(4.2,(i) and (ii)) hold when $X$ is a finite $R$-module and $x$ a
$X$--regular sequence (cf. [{\bf 5}]). However we could not prove
the converse inequality in the general case. \vspace{.2in}

\noindent {\large \bf 5. Quasi-Gorenstein ring homomorphisms}

\vspace{.2in}

In [{\bf 3}], Avramov and Foxby have defined {\it
quasi-Gorenstein} ring homomorphisms.
\para{5.1 Definition} Let $ \phi : (R,\fm)\rightarrow (S,\fn) $ be a
local ring homomorphism. Let $D$ be the dualizing complex of $R$,
then $\phi$ is said to have {\it finite Gorenstein dimension} ,if
and only if $S$ belongs to ${\cal A}(R)$; and $\phi$ is said to be
{\it quasi-Gorenstein at} $\fn$, if and only if it has finite
Gorenstein dimension and $D \otimes^{\mathbf L}_R S$ is a
dualizing complex for $S$.

They have also proved ([{\bf 3}], 7.9) that if $\phi$ is
quasi-Gorenstein at $\fn$, then an $S$-complex $X$ is in ${\cal
A}(S)$, respectively, ${\cal B}(S)$, if and only if it is in
${\cal A}(R)$, respectively, ${\cal B}(R)$.

In this appendix, using this fact, we present connections between
Gorenstein dimensions of an $S$-complex over $R$ and $S$, when
$\phi :R \rightarrow S$ is a quasi-Gorenstein ring homomorphism.

\para{5.2 Theorem} Assume that $(R,\fm)$ is a Cohen-Macaulay local
ring with a dualizing module $D$. If $ \phi : (R,\fm)\rightarrow
(S,\fn) $ is quasi-Gorenstein at $\fn$, then the following
inequalities hold for $X \in {\cal C}(S)$.
\begin{itemize}
\item[(i)] $\Gpd_R X \leq \Gpd_S X+\Gpd_R S$.
\item[(ii)] $\Gfd_R X \leq \Gfd_S X+\Gfd_R S$.
\end{itemize}

{\it Proof.} Since $S \in {\cal A}(R)$, the module $D \otimes_R S$
represents $D \otimes_R^{\mathbf{L}} S$ ( cf. [{\bf 5} ,3.4.6]).
Then $S$ is a Cohen-Macaulay local ring with a dualizing module
and hence $\Gpd_S X$ and $\Gfd_S X$ are finite if and only if
$\Gpd_R X$ and $\Gfd_R X$ are finite.

Now let $X \in {\cal A}(S)$, then using fundamental equalities we
have
$$\begin{array}{ll}
\Gpd_R X & =\sup \{-\inf (\mathbf{R} \T{Hom}_R(X,T))|T\in{\cal
I}_0 (R)\}\\
& =\sup\{-\inf(\mathbf{R} \T{Hom}_R(X
\otimes_S^{\mathbf{L}}S,T))|T\in{\cal
I}_0 (R)\}\\
& =\sup\{-\inf(\mathbf{R} \T{Hom}_S(X,\mathbf{R} \T{Hom}_R(S,T)
|T\in{\cal I}_0 (R)\}\\
& \leq
\sup\{\inf(\mathbf{R}\T{Hom}_R(S,T))-\inf(\mathbf{R}\T{Hom}_S(X,\mathbf{R}\T{Hom}_R(S,T))
|T\in{\cal I}_0
(R)\}\\
& +\sup\{-\inf(\mathbf{R}\T{Hom}_R(S,T)|T\in{\cal I}_0 (R)\}\\
& \leq \Gpd_S X+\Gpd_R S.
\end{array}$$
And
$$\begin{array}{ll}
\Gfd_R X & =\sup \{\sup (U\otimes_R^{\mathbf{L}}X)|U\in{\cal
F}_0 (R)\}\\
&
=\sup\{\sup(U\otimes_R^{\mathbf{L}}(S\otimes_S^{\mathbf{L}}X))|U\in{\cal
F}_0 (R)\}\\
&
=\sup\{\sup((U\otimes_R^{\mathbf{L}}S)\otimes_S^{\mathbf{L}}X)|U\in{\cal
F}_0 (R)\}\\
&\leq
\sup\{\sup((U\otimes_R^{\mathbf{L}}S)\otimes_S^{\mathbf{L}}X)-\sup(U\otimes_R^{\mathbf{L}}S)
|U\in{\cal F}_0 (R)\}\\
& +\sup\{\sup(U\otimes_R^{\mathbf{L}}S)
|U\in{\cal F}_0 (R)\}\\
&\leq \Gfd_SX + \Gfd_RS.
\end{array}$$
\hfill$\square$

Avramov and Foxby have proved that if $\phi :R \rightarrow S$ is
quasi-Gorenstein at $\fn$, the inequalities of the above theorem
become equalities provided that $S$ is a finite $R$-module and $X$
a finite $S$-module.\footnote { L. Avramov informed us that this
result is a consequence of a result due to Golod.} Note that when
this is the case, the Gorenstein projective and flat dimensions
are equal to the Auslander's G-dimension.

To prove the dual of the theorem, we provide some preliminaries.

For an complex $X\in {\cal C}_{(\sqsubset)}(R)$, Christensen,
Foxby and Frankild defined the {\it small restricted covariant}
Ext{\it -dimension} of $X$ as follows.

$$\begin{array}{ll}
\T{rid}_R X &\, = \sup\{-\inf(\mathbf{R}\T{Hom}_R(T,X)|T\in {\cal
P}_0^{(f)}\} \\&\, = \sup \{ -\sup U -\inf
(\mathbf{R}\T{Hom}_R(U,X)|U\in {\cal F}^{(f)}(R) \}.
\end{array}$$

They have proved that the inequality $\T{rid}_R X \leq \id_R X$
always holds, and it becomes equality if $\T{cmd} R \leq 1$ and
$\id_R X < \infty$.({\bf cf}.[{\bf 7}])

\para{5.3 Theorem} Let $M$ be an $R$-module. The following inequality
always holds.
$$\T{rid}_ R M \leq \Gid_R M .$$

The equality holds if $\T{cmd} R \leq 1$ and $\Gid_R M < \infty$.

{\it Proof.} If $\Gid_R M$ is not finite, then the claim is
obvious. Now let $n=\Gid_R M $ be finite. We prove the theorem by
induction on $n$.

If $n=0$, then $M$ is a Gorenstein injective module and then there
exists an exact complex
$$ I=\cdots \rightarrow I_1 \rightarrow I_0 \rightarrow M
\rightarrow 0 $$ with $I_j$ s injective modules .

If $T$ is a finite $R$-module of finite projective dimension $t$,
then since $\T{Ext}_R^i(T,M)= \T{Ext}_R^{i+t} (T,\T{Z_t^I})$, we
have $\T{Ext}_R^i(T,M)=0$ for all $i>0$. Therefore $\T{rid}_R M =
0$.

If $n>0$, then by [{\bf 17}, 2.45], there exists a Gorenstein
injective $R$-module $G$ and an $R$-module $C$ with $\id_R C
=\Gid_R C = n-1$, such that the following sequence is exact.
$$ 0 \rightarrow M \rightarrow G \rightarrow C \rightarrow 0$$

For any finite $R$-module $T$ of finite projective dimension, we
have the exact sequence
$$0=\T{Ext}_R^{i-1}(T,C) \rightarrow \T{Ext}_R^i(T,M) \rightarrow
\T{Ext}_R^i(T,G)=0$$ for $i>n$.

Therefore $\T{rid}_RM \leq \Gid_R M$.

Now let $ \T{cmd}R \leq 1 $. To prove the inverse inequality, note
that since $G$ is a Gorenstein injective module, there exists an
injective $R$-module $E$ and an exact sequence
$$0 \rightarrow K \rightarrow E \rightarrow G \rightarrow 0 $$
with $K$ a Gorenstein injective $R$-module, too.

On the other hand we have the isomorphisms $C \cong G/M$ and $G
\cong E/K$. Then there exists a submodule $L$ of $E$ such that $K
\subseteq L$ and $M \cong L/K$, and then $C \cong E/L$.

From the following exact sequence, we get $\id_R L \leq \id_R
C=n$.
$$0 \rightarrow L \rightarrow E \rightarrow C \rightarrow 0 $$

Now consider the exact sequence
$$0 \rightarrow K \rightarrow L \rightarrow M \rightarrow 0 .$$
Since $\Gid_R M =n$, there exists an injective $R$-module $J$ with
$\T{Ext}_R^n(J,M)\neq 0$. Hence from the exact sequence
$$\T{Ext}_R^n(J,L) \rightarrow \T{Ext}_R^n(J,M) \rightarrow
\T{Ext}_R^{n+1}(J,K)=0,$$ we get $\T{Ext}_R^n(J,L)\neq 0$ and then
$\id_R L \geq n$.

Therefore $\T{rid}_R L=\id_R L=n$, and then there exists a finite
$R$-module $Q$ of finite projective dimension such that
$\T{Ext}_R^n(Q,L)\neq 0$.

Finally the exactness of
$$0=\T{Ext}_R^n(Q,K) \rightarrow \T{Ext}_R^n(Q,L) \rightarrow
\T{Ext}_R^n(Q,M)$$ implies that $\T{Ext}_R^n(Q,M)\neq 0$ and then
$\T{rid}_R M \geq n$. \hfill$\square$

\para{5.4 Theorem } Assume that $(R,\fm)$ is a Cohen-Macaulay local
ring with a dualizing module $D$. If $ \phi : (R,\fm)\rightarrow
(S,\fn) $ is quasi-Gorenstein at $\fn$, then the following
inequalities hold for an $S$-module $M$.
$$\Gid_R M \leq \Gid_S M + \Gfd_R S.$$
Furthermore, the equality holds if $S$ is a finite $R$-module and
$M$ a finite $S$-module.

{\it Proof.} Since $\phi$ is quasi-Gorenstein, $\Gid_R M$ is
finite, if and only if $\Gid_S M$ is finite. Let $\Gid_S M <
\infty$, then
$$\begin{array}{ll}
\Gid_R M & =\sup \{-\inf(\mathbf{R}\T{Hom}_R(T,M)|T\in{\cal
P}_0^{(f)} (R)\}\\
&=\sup\{-\inf(\mathbf{R}\T{Hom}_R(T,\mathbf{R}\T{Hom}_S(S,M)))|T\in{\cal
P}_0^{(f)} (R)\}\\
&=\sup\{-\inf(\mathbf{R}\T{Hom}_S(T\otimes_R^{\mathbf{L}}S,M))|T\in{\cal
P}_0^{(f)} (R)\}\\
&\leq
\sup\{-\sup(T\otimes_R^{\mathbf{L}}S)-\inf(\mathbf{R}\T{Hom}_S(T
\otimes_R^{\mathbf{L}}S,M))|T\in{\cal P}_0^{(f)} (R)\}\\
&+\sup\{\sup(T\otimes_R^{\mathbf{L}}S)|T\in{\cal P}_0^{(f)} (R)\}\\
& \leq \Gid_S M + \Gfd_R S.
\end{array}$$

Now let $S$ be a finite $R$-module and $M$ a finite $S$-module. By
[{\bf 5}, 6.2.15], the following equalities hold.
$$\Gid_R M=\T{depth}R \ \ \T{and} \ \ \Gid_S M=\T{depth}S.$$

And then the requested equality follows by  Auslander-Bridger
formula. \hfill$\square$

\vspace{.3in}





\baselineskip=18pt

\noindent {\large \bf References}

\noindent [1]
\begin{minipage}[t]{5.5in}
M.\ Auslander, {\em Anneaux de Gorenstein et torsion en
alg\`{e}bre commutative}, S\'{e}minaire d'alg\`{e}bre commutative
1966/67, notes by M.\ Mangeney, C.\ Peskine and L. Szpiro,
\'{E}cole Normale Sup\'{e}rieure de Jeunes Filles, Paris, 1967.
\end{minipage}
\vspace{.1in}

\noindent [2]
\begin{minipage}[t]{5.5in}
L.L.\ Avramov and H.-B.\ Foxby, {\em Homological dimensions of
unbounded complexes}, J.\ Pure Appl.\ Algebra {\bf 71} (1991),
129-155.
\end{minipage}
\vspace{.1in}

\noindent [3]
\begin{minipage}[t]{5.5in}
L.L.\ Avramov and H.-B.\ Foxby, {\em Ring homomorphisms and finite
Gorenstein dimension}, Proc.\ London\ Math.\ Soc.\ (3) {\bf 75}
(1997), no. 2, 241--270.
\end{minipage}
\vspace{.1in}

\noindent [4]
\begin{minipage}[t]{5.5in}
W.\ Burns and J.\ Herzog, {\em Cohen--Macaulay rings}, Cambridge
Studies in Advanced Mathematics, {\bf 39}, Cambridge University
press, Cambridge, 1993.
\end{minipage}
\vspace{.1in}

\noindent [5]
\begin{minipage}[t]{5.5in}
L.W.\ Christensen, {\em Gorenstein dimensions}, Lecture notes in
Mathematics 1747, Springer-Verlag, Berlin, 2000.
\end{minipage}
\vspace{.1in}

\noindent [6]
\begin{minipage}[t]{5.5in}
L.W.\ Christensen, {\em Semi-dualizing complexes and their
Auslander categories}, Trans.\ Amer.\ Math.\ Soc.\ {\bf 353}
(2001), no.\ 5, 1839-1883.
\end{minipage}
\vspace{.1in}

\noindent [7]
\begin{minipage}[t]{5.5in}
L.W.\ Christensen, H.-B.\ Foxby and A.\ Frankild, {\em Functorial
dimensions and Cohen-Macaulayness of rings}, To appear in J.\
Algebra.
\end{minipage}
\vspace{.1in}

\noindent [8]
\begin{minipage}[t]{5.5in}
E.E.\ Enochs and O.M.G.\ Jenda, {\em On Gorenstein injective
modules}, Comm.\ Algebra {\bf 21} (1993), no.\ 10, 3489-3501.
\end{minipage}
\vspace{.1in}

\noindent [9]
\begin{minipage}[t]{5.5in}
E.E.\ Enochs and O.M.G.\ Jenda, {\em Gorenstein injective and
projective modules}, Math.\ Z.\ {\bf 220} (1995), no.\ 4, 611-633.
\end{minipage}
\vspace{.1in}

\noindent [10]
\begin{minipage}[t]{5.5in}
E.E.\ Enochs and O.M.G.\ Jenda, {\em Gorenstein injective and flat
dimensions}, Math.\ Japan {\bf 44} (1996), no.\ 2, 261-268.
\end{minipage}
\vspace{.1in}

\noindent [11]
\begin{minipage}[t]{5.5in}
E.E.\ Enochs and O.M.G.\ Jenda, {\em Gorenstein injective modules
over Gorenstein rings}, Comm.\ Algebra {\bf 26} (1998), no.\ 11,
3489-3496.
\end{minipage}
\vspace{.1in}

\noindent [12]
\begin{minipage}[t]{5.5in}
E.E.\ Enochs, O.M.G.\ Jenda and B.\ Torrecilles, {\em Gorenstein
flat modules}, Nanjing Daxue Xuebao Shuxue Bannian Kan {\bf 10}
(1993), no.\ 1, 1-9.
\end{minipage}
\vspace{.1in}

\noindent [13]
\begin{minipage}[t]{5.5in}
E.E.\ Enochs, O.M.G.\ Jenda and Jinzhong Xu, {\em Foxby duality
and Gorenstein injective and projective modules}, Trans.\ Amer.\
math.\ Soc.\ {\bf 348} (1996), no.\ 8, 3223-3234.
\end{minipage}
\vspace{.1in}

\noindent [14]
\begin{minipage}[t]{5.5in}
H.-B.\ Foxby, {\em Hyperhomological algebra and commutative
rings}, notes in preparation.
\end{minipage}
\vspace{.1in}

\noindent [15]
\begin{minipage}[t]{5.5in}
H.-B.\ Foxby, {\em Gorenstein dimension over Cohen--Macaulay
rings}, Proceeding of International Conference on Commutative
Algebra (W.\ Bruns.\ ed.), Universit\"{a}t Onsabr\"{u}ck, 1994.
\end{minipage}
\vspace{.1in}

\noindent [16]
\begin{minipage}[t]{5.5in}
H.\ Holm, {\em Gorenstein homological dimensions}, To appear.
\end{minipage}
\vspace{.1in}

\noindent [17]
\begin{minipage}[t]{5.5in}
T.\ Ishikawa, {\em On injective modules and flat modules}, J.\
Math.\ Soc.\ Japan {\bf 17} (1965), 291-296.
\end{minipage}
\vspace{.1in}

\noindent [18]
\begin{minipage}[t]{5.5in}
S.\ Yassemi, {\em On flat and injective dimensions}, Ital.\ J.\ of
Pure Appl.\ Math.\ No.\ 6 (1999), 33-41.
\end{minipage}
\vspace{.1in}

\noindent [19]
\begin{minipage}[t]{5.5in}
S. Yassemi, {\em Homomorphism locally of finite injective
dimension}, New Zealand J.\ Math.\ {\bf 26} (2000), no.\ 1,
97-102.
\end{minipage}
\vspace{.1in}

\end{document}